\documentclass[12pt]{amsart}

\usepackage{amsmath,xspace,amssymb,mathrsfs}
\usepackage{color}

\input xy
\xyoption{all}
\xyoption{2cell}
\UseAllTwocells
\CompileMatrices

\newcommand{\Spec}{\operatorname{Spec}}
\renewcommand{\phi}{\varphi}

\newcommand{\Ima}{\operatorname{Im}}

\newcommand{\Max}{\operatorname{Max}}

\newcommand{\Sp}{\operatorname{Sq}}

\newtheorem{proposition}{Proposition}[section]
\newtheorem{lemma}[proposition]{Lemma} 
\newtheorem{corollary}[proposition]{Corollary}
\newtheorem{theorem}[proposition]{Theorem}

\newtheorem{prop-def}[proposition]{Proposition and definition}

\theoremstyle{definition}
\newtheorem{definition}[proposition]{Definition}

\newtheorem{remark}[proposition]{Remark}

\begin{document}

\title{Quasi-prime ideals}

\author[A. Tarizadeh and M.Aghajani]{Abolfazl Tarizadeh and Mohsen Aghajani}
\address{ Department of Mathematics, Faculty of Basic Sciences, University of Maragheh \\
P. O. Box 55136-553, Maragheh, Iran.
 }
\email{ebulfez1978@gmail.com, aghajani14@gmail.com}


\begin{abstract} In this paper, the new concept of quasi-prime ideal is introduced which at the same time generalizes the ``prime ideal'' and ``primary ideal'' notions. Then a natural topology on the set of quasi-prime ideals of a ring is introduced which generalizes the Zariski topology. The basic properties of the quasi-prime spectrum are studied and several interesting results are obtained. Specially, it is proved that if the Grothendieck t-functor is applied on the quasi-prime spectrum then the prime spectrum is deduced. It is also shown that there are the cases that the prime spectrum and quasi-prime spectrum do not behave similarly.

\end{abstract}

\maketitle

\section{Introduction}

In this paper, we first generalize naturally the notion of prime ideal. In fact, a proper ideal $\mathfrak{q}$ of a commutative ring $A$ is said to be a \emph{quasi-prime} ideal of $A$ if $fg\in\mathfrak{q}$ for some $f,g\in A$ then either $f\in\sqrt{\mathfrak{q}}$ or $g\in\sqrt{\mathfrak{q}}$. The set of quasi-prime ideals of $A$ is denoted by $\Sp A$. Then we equip this set with a natural topology whose basis opens are of the form $U_{f}=\{\mathfrak{q}\in\Sp A: \mathfrak{q}\cap S_{f}=\emptyset\}$ where $f\in A$ and $S_{f}=\{1,f,f^{2},...\}$. The space $\Sp A$ is called the quasi-prime spectrum of $A$. It is shown that the quasi-prime spectrum satisfies in all of the conditions of being a spectral space except the uniqueness of generic point, see Theorem \ref{Lemma II} and Proposition \ref{Proposition II}. This topology is a natural generalization of the Zariski topology, i.e. the prime spectrum $\Spec A$ is a dense subspace of $\Sp A$. In many cases, this topology behaves completely different from the Zariski topology. For instance, $\Sp\mathbb{Z}$ has no closed point. The basic properties of the quasi-prime spectrum are studied. Indeed, Theorems \ref{Lemma II}, \ref{Lemma III}, \ref{Theorem II}, Proposition \ref{Proposition IV}, Theorem \ref{Theorem III}, Propositions \ref{Theorem I} and \ref{Proposition II} are the main results on this topology. In Theorem \ref{Theorem III}, the connected components of the quasi-prime spectrum are characterized. In Theorem \ref{Theorem II}, it is shown that the prime spectrum can be canonically recovered from the quasi-prime spectrum by applying the Grothendieck t-functor. \\

\section{Quasi-primes}

\begin{definition} A proper ideal $\mathfrak{q}$ of a ring $A$ is called a quasi-prime ideal of $A$ if $fg\in\mathfrak{q}$ for some $f,g\in A$ then either $f\in\sqrt{\mathfrak{q}}$ or $g\in\sqrt{\mathfrak{q}}$. \\
\end{definition}

\begin{lemma} An ideal $\mathfrak{q}$ of a ring $A$ is a quasi-prime ideal of $A$ if and only if $\sqrt{\mathfrak{q}}$ is a prime ideal of $A$.\\
\end{lemma}

{\bf Proof.} Easy. $\Box$ \\

Recall that a proper ideal $\mathfrak{q}$ of a ring $A$ is called a primary ideal of $A$ if there exist $f,g\in A$ such that $fg\in\mathfrak{q}$ and $f\notin\mathfrak{q}$ then $g\in\sqrt{\mathfrak{q}}$. It is equivalent to the statement that if there exist $f,g\in A$ such that $fg\in\mathfrak{q}$ and $f\notin\sqrt{\mathfrak{q}}$ then $g\in\mathfrak{q}$. Note that if $\mathfrak{q}$ is a primary ideal of $A$, $fg\in\mathfrak{q}$, $f\notin\mathfrak{q}$ and $g\notin\mathfrak{q}$ then by definition both of $f$ and $g$ must be in $\sqrt{\mathfrak{q}}$. Every primary ideal is a quasi-prime ideal but the converse is not necessarily true. As a specific example, let $A$ be the polynomial ring $k[x,y,z]$ modulo $I=(xy-z^{2})$ where $k$ is a domain. Then $\mathfrak{p}=(x+I, z+I)$ is a prime ideal of $A$ since $A/\mathfrak{p}\simeq k[y]$, but $\mathfrak{q}=\mathfrak{p}^{2}$ is a quasi-prime ideal which is not a primary ideal. Because $(x+I)(y+I)=(z+I)^{2}\in\mathfrak{q}$ but $x+I\notin\mathfrak{q}$ and $y+I\notin\sqrt{\mathfrak{q}}$. \\

If $\mathfrak{p}$ is a prime ideal of a ring $A$ then $\mathfrak{p}^{n}$ is a quasi-prime ideal of $A$ for all $n\geq1$. But there are also non-primary quasi-prime ideals which are not as the power of a prime ideal. As an example,
let $A$ be the polynomial ring $k[x,y,z,t]$ modulo $I=(xy-z^{2})$ where $k$ is a domain. Then $\mathfrak{p}=(z+I,t+I)$ is a prime of $A$ since $A/\mathfrak{p}\simeq k[x,y]$, but $\mathfrak{q}=(z+I, t^{2}+I)$ is a non-primary quasi-prime of $A$ which is not also as a power of a prime ideal. \\

The set of quasi-prime ideals of $A$ is denoted by $\Sp A$. For each $f\in A$ we define $U_{f}=\{\mathfrak{q}\in\Sp A: \mathfrak{q}\cap S_{f}=\emptyset\}$ where $S_{f}=\{f^{i}:i\geq0\}=\{1,f,f^{2},...\}$. Then clearly $U_{1}=\Sp A$ and $U_{f}\cap U_{g}=U_{fg}$ for all $f,g\in A$. Thus there exists a (unique) topology over $\Sp A$ such that the basis opens are precisely of the form $U_{f}$ where $f\in A$. We call $\Sp A$ the quasi-prime spectrum (or, the space of quasi-primes) of $A$. Clearly $\Spec A$ is a subspace of $\Sp A$ since $D(f)=U_{f}\cap\Spec A$ for all $f\in A$. Note that $U_{f}=\{\mathfrak{q}\in\Sp A: \sqrt{\mathfrak{q}}\cap S_{f}=\emptyset\}$ for all $f\in A$. It follows that $\Spec A$ is a dense subspace of $\Sp A$. \\

\begin{theorem}\label{Lemma II} Every basis open $U_{f}$ is quasi-compact. In particular, $\Sp A$ is quasi-compact. \\
\end{theorem}

{\bf Proof.} It suffices to show that every open covering of $U_{f}$ by the basis opens has a finite refinement. Hence let $U_{f}=\bigcup\limits_{i\in I}U_{g_{i}}$ where $g_{i}\in A$ for all $i$. It follows that $f\in\sqrt{(g_{i}: i\in I)}$. Thus the exists a finite subset $J$ of $I$ such that $f\in\sqrt{(g_{i}: i\in J)}$. We show that $U_{f}=\bigcup\limits_{i\in J}U_{g_{i}}$. If $\mathfrak{q}\in U_{f}$ then there exists some $i\in J$ such that $g_{i}\notin\sqrt{\mathfrak{q}}$. It follows that $\mathfrak{q}\in U_{g_{i}}$. $\Box$ \\

The proof of Theorem \ref{Lemma II} also shows that $D(f)$ is quasi-compact for all $f\in A$. \\

\begin{corollary} If $X$ is a subspace of $\Sp A$ such that $\Spec A\subseteq X$ then $X$ is quasi-compact. In particular, the primary spectrum (space of primary ideals) of $A$ is quasi-compact. \\
\end{corollary}

{\bf Proof.} If $U$ is an open of $\Sp A$ such that $X\subseteq U$ then $U=\Sp A$. Because if $\mathfrak{q}$ is a quasi-prime ideal of $A$ then there exists some $f\in A$ such that $\sqrt{\mathfrak{q}}\in U_{f}\subseteq U$. It follows that $\mathfrak{q}\in U_{f}$. Thus by Theorem \ref{Lemma II}, $X$ is quasi-compact.
$\Box$ \\

\begin{theorem}\label{Lemma III} If $\mathfrak{q}$ is a quasi-prime ideal of $A$ then $\overline{\{\mathfrak{q}\}}=\{\mathfrak{p}\in\Sp A: \mathfrak{q}\subseteq\sqrt{\mathfrak{p}}\}$. \\
\end{theorem}

{\bf Proof.} Let $\mathfrak{p}\in\overline{\{\mathfrak{q}\}}$ and $f\in\mathfrak{q}$. If $f\notin\sqrt{\mathfrak{p}}$ then $\mathfrak{p}\in U_{f}$. It follows that $\mathfrak{q}\in U_{f}$, a contradiction. Conversely, assume that $\mathfrak{q}\subseteq\sqrt{\mathfrak{p}}$. If $\mathfrak{p}\notin\overline{\{\mathfrak{q}\}}$ then there exists some $f\in A$ such that $\mathfrak{p}\in U_{f}$ but $\mathfrak{q}\notin U_{f}$. Hence there exists a natural number $n\geq1$ such that $f^{n}\in\mathfrak{q}$. It follows that $f\in\sqrt{\mathfrak{p}}$. But this is contradiction since $\mathfrak{p}\cap S_{f}=\emptyset$. $\Box$ \\

Unlike the prime spectrum, a maximal ideal is not necessarily a closed point of the quasi-prime spectrum. As a specific example, if $p$ is a prime number then $\{p^{n}\mathbb{Z}: n\geq1\}$ is the closure of $\{p\mathbb{Z}\}$ in $\Sp\mathbb{Z}$. In fact, $\Sp\mathbb{Z}$ has no closed point. \\

In the space $\Sp A$ generic points are not unique:\\

\begin{corollary}\label{Corollary II} If $\mathfrak{q}$ is a quasi-prime ideal of $A$ then $\overline{\{\mathfrak{q}\}}=\overline{\{\sqrt{\mathfrak{q}}\}}$. $\Box$ \\
\end{corollary}

\begin{corollary}\label{Corollary I} Let $E$ be a closed subset of $\Sp A$. If $\mathfrak{q}\in E$ then $\sqrt{\mathfrak{q}}\in E$. $\Box$ \\
\end{corollary}

It is easy to see that the closed subsets of $\Sp A$ are precisely of the form $\mathcal{V}(I)=\{\mathfrak{q}\in\Sp A: I\subseteq\sqrt{\mathfrak{q}}\}$ where $I$ is an ideal of $A$. Clearly $\mathcal{V}(I)\cap\Spec A=V(I)$.
If $f\in A$ then $\mathcal{V}(f)=\{\mathfrak{q}\in\Sp A: f\in\sqrt{\mathfrak{q}}\}$. If $\mathfrak{p}$ is a quasi-prime of $A$ then by Theorem \ref{Lemma III}, $\mathcal{V}(\mathfrak{p})=\overline{\{\mathfrak{p}\}}$. \\

\section{Connected components}

A subspace $Y$ of a topological space $X$ is called a retraction of $X$ if there exists a continuous map $\gamma:X\rightarrow Y$ such that $\gamma(y)=y$ for all $y\in Y$. Such a map $\gamma$ is called a retraction map. \\

\begin{lemma}\label{Lemma I} The prime spectrum is a retraction of quasi-prime spectrum.\\
\end{lemma}

{\bf Proof.} The map $\gamma:\Sp A\rightarrow\Spec A$ given by $\mathfrak{q}\rightsquigarrow\sqrt{\mathfrak{q}}$ is continuous since $\gamma^{-1}\big(D(f)\big)=U_{f}$ for all $f\in A$. $\Box$ \\

The map $\gamma$ in Lemma \ref{Lemma I} is an open map and $\gamma^{-1}\big(V(I)\big)=\mathcal{V}(I)$ for all ideals $I$ of $A$. \\

\begin{remark}\label{Remark II} There is a fundamental result due to Grothendieck which states that the map $f\rightsquigarrow D(f)$ is a bijection from the set of idempotents of $A$ onto the set of clopen (both open and closed) subsets of $\Spec A$, see \cite[Tag 00EE]{Johan}. Under the light of this result we obtain that: \\
\end{remark}

\begin{proposition}\label{Proposition IV} The map $f\rightsquigarrow U_{f}$ is a bijection from the set of idempotents of $A$ onto the set of clopens of $\Sp A$. \\
\end{proposition}

{\bf Proof.} By Remark \ref{Remark II}, it suffices to show that the map $U\rightsquigarrow\gamma^{-1}(U)$ is a bijection from the set of clopens of $\Spec A$ onto the set of clopens of $\Sp A$, for $\gamma$ see Lemma \ref{Lemma I}. Assume that $\gamma^{-1}(U)=\gamma^{-1}(V)$. If $\mathfrak{p}\in U$ then $\mathfrak{p}\in\gamma^{-1}(U)$ and so $\gamma(\mathfrak{p})=\mathfrak{p}\in V$. Hence $U=V$. It remains to show that this map is surjective. If $U$ is a clopen of $\Sp A$ then $U\cap\Spec A$ is a clopen of $\Spec A$. We have $U=\gamma^{-1}(U\cap\Spec A)$. Because if $\mathfrak{q}\in U$ then by Corollary \ref{Corollary I}, $\gamma(\mathfrak{q})=\sqrt{\mathfrak{q}}\in U\cap\Spec A$. Conversely, if $\mathfrak{q}\in\gamma^{-1}(U\cap\Spec A)$ then there exists some $f\in A$ such that $\sqrt{\mathfrak{q}}\in U_{f}\subseteq U$. It follows that $\mathfrak{q}\in U$. $\Box$ \\

\begin{corollary}\label{Corollary IV} The space $\Sp A$ is connected if and only if $A$ has no nontrivial idempotents. $\Box$ \\
\end{corollary}

\begin{proposition}\label{Proposition III} If $\phi:A\rightarrow B$ is a morphism of rings then the induced map $\phi^{\ast}:\Sp B\rightarrow \Sp A$ given by $\mathfrak{q}\rightsquigarrow\phi^{-1}(\mathfrak{q})$ is continuous. \\
\end{proposition}

{\bf Proof.} If $\mathfrak{q}$ is a quasi-prime ideal of $B$ then $\phi^{-1}(\mathfrak{q})$ is a quasi-prime ideal of $A$ because $\phi^{-1}(\sqrt{\mathfrak{q}})=\sqrt{\phi^{-1}(\mathfrak{q})}$. Hence $\phi^{\ast}$ is well-defined. It is continuous since
$(\phi^{\ast})^{-1}(U_{f})=U_{\phi(f)}$ for all $f\in A$. $\Box$ \\

Note that if $S$ is a multiplicative subset of $A$ then the map $\pi^{\ast}:\Sp (S^{-1}A)\rightarrow\Sp A$ induced by the canonical ring map $A\rightarrow S^{-1}A$ is not injective and $\Ima\pi^{\ast}\subseteq\{\mathfrak{q}\in\Sp A: \mathfrak{q}\cap S=\emptyset\}$. Specially $\Ima\pi^{\ast}\subseteq U_{f}$ where $\pi:A\rightarrow A_{f}$ is the canonical map. \\

\begin{lemma}\label{Lemma IV} If $I$ is an ideal of a ring $A$ then the map $\pi:\Sp A/I\rightarrow \Sp A$ induced by the canonical ring map $A\rightarrow A/I$ is injective and $\Ima\pi=\mathcal{V}(I)$. \\
\end{lemma}

{\bf Proof.} Easy. $\Box$ \\

An ideal of $A$ is said to be a regular ideal of $A$ if it is generated by a subset of idempotent elements of $A$. Each maximal element of the set of proper regular ideals of $A$ (ordered by inclusion) is called a max-regular ideal of $A$. By the Zorn's Lemma, every proper regular ideal of $A$ is contained in a max-regular ideal of $A$. It is well known that a regular ideal $M$ is a max-regular ideal of $A$ if and only if $A/M$ has no nontrivial idempotents, see \cite[Lemma 3.19]{Abolfazl}. It is also well known that the connected components of $\Spec A$ are precisely of the form $V(M)$ where $M$ is a max-regular ideal of $A$, see \cite[Theorem 3.17]{Abolfazl}. We have then the following result. \\

\begin{theorem}\label{Theorem III} The connected components of $\Sp A$ are precisely of the form $\mathcal{V}(M)$ where $M$ is a max-regular ideal of $A$. \\
\end{theorem}

{\bf Proof.} If $N$ is a max-regular ideal of $A$. Then by Corollary \ref{Corollary IV}, $\Sp A/N$ is connected. Thus by Proposition \ref{Proposition III} and Lemma \ref{Lemma IV}, $\mathcal{V}(N)$ is connected. Now let $C$ be a connected component of $\Sp A$. Then $\gamma(C)$ is contained in a connected component of $\Spec A$, for $\gamma$ see Lemma \ref{Lemma I}. Thus there exists a max-regular ideal $M$ of $A$ such that $\gamma(C)\subseteq V(M)$. It follows that $C\subseteq\mathcal{V}(M)$. Thus $C=\mathcal{V}(M)$ since $\mathcal{V}(M)$ is connected. Conversely, let $N$ be a max-regular ideal of $A$. Then there exists a connected component $C$ of $\Sp A$ such that $\mathcal{V}(N)\subseteq C$. We observed that there exists a max-regular ideal $M$ of $A$ such that $C=\mathcal{V}(M)$. It follows that $\sqrt{M}\subseteq\sqrt{N}$. Thus $M\subseteq N$ since $M$ is a regular ideal. This implies that $M=N$ because $M$ is max-regular and $N$ is regular ideal. $\Box$ \\

\section{t-functor}

There exists a covariant functor due to Grothendieck from the category of topological spaces to itself. It is called the t-functor. This functor has geometric applications and builds a bridge between the classical algebraic geometry and modern algebraic geometry. In what follows we shall introduce this functor. If $X$ is a topological space then the points of $t(X)$ are the irreducible and closed subsets of $X$. Recall that a topological space is said to be an irreducible space if it is non-empty and can not be written as the union of two proper closed subsets. The closed subsets of $t(X)$ are precisely of the form $t(Y)$ where $Y$ is a closed subset of $X$. If $f:X\rightarrow X'$ is a continuous map of topological spaces then the function $t(f):t(X)\rightarrow t(X')$ given by $Z\rightsquigarrow\overline{f(Z)}$ is well-defined and continuous. There exists also a canonical continuous map $X\rightarrow t(X)$ defined by $x\rightsquigarrow \overline{\{x\}}$. We have then the following result. \\

\begin{theorem}\label{Theorem II} The space $t(\Sp A)$ is canonically homeomorphic to $\Spec A$. \\
\end{theorem}

{\bf Proof.} Let $Z$ be an irreducible and closed subset of $\Sp A$. By Proposition \ref{Proposition II}, there exists a quasi-prime $\mathfrak{q}$ of $A$ such that $Z=\overline{\{\mathfrak{q}\}}$. We then define $\phi:t(\Sp A)\rightarrow\Spec A$ as $Z\rightsquigarrow\sqrt{\mathfrak{q}}$. We show that it is a homeomorphism. The map $\phi$ is injective, see Corollary \ref{Corollary II}. If $\mathfrak{p}$ is a prime ideal of $A$ then $Z:=\mathcal{V}(\mathfrak{p})=\overline{\{\mathfrak{p}\}}$ is an irreducible and closed subset of $\Sp A$ and $\phi(Z)=\mathfrak{p}$. The map $\phi$ is continuous because $\phi^{-1}\big(D(f)\big)=t(\Sp A)\setminus t\big(\mathcal{V}(f)\big)$ for all $f\in A$. It remains to show that $\phi$ is a closed map. Let $Y=\mathcal{V}(I)$ be a closed subset of $\Sp A$ where $I$ is an ideal of $A$. Then $\phi\big(t(Y)\big)=V(I)$. $\Box$ \\

\section{Normality}

A topological space is called a normal space if every two disjoint closed subsets admit disjoint open neighborhoods. Clearly a closed subspace of a normal space is a normal space, but an arbitrary subspace is not necessarily a normal space. \\

\begin{proposition}\label{Theorem I} Let $A$ be a ring. Then $\Sp A$ is a normal space if and only if $\Spec A$ is a normal space. \\
\end{proposition}

{\bf Proof.} Let $\Sp A$ be a normal space. Let $E=V(I)$ and $F=V(J)$ be two disjoint closed subsets of $\Spec A$ where $I$ and $J$ are ideals of $A$. It follows that $I+J=A$. Thus $\mathcal{V}(I)\cap\mathcal{V}(J)=\emptyset$.
Hence there are disjoint opens $U$ and $V$ in $\Sp A$ such that $\mathcal{V}(I)\subseteq U$ and $\mathcal{V}(J)\subseteq V$. It follows that $E\subseteq U\cap\Spec A$ and $F\subseteq V\cap\Spec A$. Conversely, let $\Spec A$ be a normal space. Let $\mathcal{V}(I)$ and $\mathcal{V}(J)$ be two disjoint closed subsets of $\Sp A$. It follows that $V(I)\cap V(J)=\emptyset$. Thus there are disjoint opens $U$ and $V$ in $\Spec A$ such that $V(I)\subseteq U$ and $V(J)\subseteq V$. It follows that $\mathcal{V}(I)\subseteq\gamma^{-1}(U)$ and $\mathcal{V}(I)\subseteq\gamma^{-1}(U)$. Hence $\Sp A$ is a normal space. $\Box$ \\

It is well known that $\Spec A$ is a normal space if and only if $A$ is a pm-ring, see \cite[Theorem 1.2]{Marco-Orsatti} or \cite[Theorem 4.3]{Abolfazl-Mohsen}. \\

\begin{proposition}\label{Proposition I} If $\Spec A$ is a normal space then $\Max A$ is a normal space. \\
\end{proposition}

{\bf Proof.} It is well known that $\Max A$ is quasi-compact. It suffices to show that $\Max A$ is Hausdorff because it is well known that every compact (quasi-compact and Hausdorff) space is a normal space. Thus let $\mathfrak{m}$ and $\mathfrak{m}'$ be two distinct maximal ideals of $A$. The closed points of $\Spec A$ are precisely the maximal ideals. Therefore by the hypothesis, there exist disjoint opens $U$ and $V$ in $\Spec A$ such that $\mathfrak{m}\in U$ and $\mathfrak{m}'\in V$. It follows that $(U\cap\Max A)\cap(V\cap\Max A)=\emptyset$. $\Box$ \\

The converse of Proposition \ref{Proposition I} is not necessarily true, see \cite[Remark 4.8]{Abolfazl-Mohsen}. This Remark, as stated there, also shows that the main result of \cite{Simmons} is not true. \\

\section{Spectrality and Hausdorfness}

\begin{proposition}\label{Proposition II} Every irreducible and closed subset of $\Sp A$ has a generic point. \\
\end{proposition}

{\bf Proof.} Let $Z$ be an irreducible and closed subset of $\Sp A$. There exists an ideal $J$ of $A$ such that $Z=\mathcal{V}(J)$. By Theorem \ref{Lemma III}, it suffices to show that $J$ is a quasi-prime of $A$. Clearly $J\neq A$ since $Z$ is non-empty. Let $f,g\in A$ such that $fg\in J$. We have $Z=\big(\mathcal{V}(f)\cap Z\big)\cup\big(\mathcal{V}(g)\cap Z\big)$. \\
It follows that either $Z\subseteq\mathcal{V}(f)$ or $Z\subseteq\mathcal{V}(g)$. Thus either $f\in\sqrt{J}$ or $g\in\sqrt{J}$. $\Box$ \\

The irreducible components of $\Sp A$ are precisely of the form $\mathcal{V}(\mathfrak{p})$ where $\mathfrak{p}$ is a minimal prime of $A$.\\

Recall that a topological space $X$ is called a spectral space if it is quasi-compact, its topology has a basis consisting of quasi-compact opens such that every finite intersection of these basis opens is again quasi-compact and every irreducible and closed subset of $X$ has a unique generic point. In this definition, the uniqueness of generic point is a crucial point. For instance,
the space $\Sp A$, by Theorem \ref{Lemma II} and Proposition \ref{Proposition II}, satisfies in all of the conditions of a spectral space except the uniqueness of generic point. This leads us to the following result: \\

\begin{corollary} $\Sp A$ is a spectral space if and only if $\Sp A=\Spec A$. $\Box$ \\
\end{corollary}

\begin{remark}\label{Remark I} It is well known that for a ring $A$ then $\Spec A$ is Hausdorff if and only if every prime ideal of $A$ is maximal. \\
\end{remark}

\begin{proposition} $\Sp A$ is a Hausdorff space if and only if $\Sp A=\Max A$. \\
\end{proposition}

{\bf Proof.} Let $\Sp A$ be a Hausdorff space and $\mathfrak{q}$ a quasi-prime of $A$. There exists a maximal ideal $\mathfrak{m}$ of $A$ such that $\mathfrak{q}\subseteq\mathfrak{m}$. It follows that $\mathfrak{m}\in\mathcal{V}(\mathfrak{q})=\{\mathfrak{q}\}$ and so $\mathfrak{q}=\mathfrak{m}\in\Max A$. The converse implies from Remark \ref{Remark I}. $\Box$ \\

\begin{corollary}\label{Corollary III} If $\Sp A$ is a Hausdorff space then $\Spec A$ is Hausdorff. $\Box$ \\
\end{corollary}

The converse of Corollary \ref{Corollary III} does not hold. As a specific example, if $A=\mathbb{Z}/8\mathbb{Z}$ then
$\Spec A=\{\mathfrak{m}\}$ but $\Sp A=\{0,\mathfrak{m},\mathfrak{m}^{2}\}$ where $\mathfrak{m}=2\mathbb{Z}/8\mathbb{Z}$. \\

\begin{corollary} Let $A$ be a local ring such that the maximal ideal $\mathfrak{m}$ is a finitely generated ideal. Then the following are equivalent. \\
$\mathbf{(i)}$ $A$ is a field. \\
$\mathbf{(ii)}$ $\Sp A$ is Hausdorff. \\
$\mathbf{(iii)}$ $\Sp A$ is a spectral space. \\
$\mathbf{(iv)}$ $\Sp A$ has a closed point. \\
\end{corollary}

{\bf Proof.} If one of the statements $(ii),(iii)$ and $(iv)$ hold then we have $\mathfrak{m}=\mathfrak{m}^{2}$. Thus by the Nakayama lemma, $\mathfrak{m}=0$. $\Box$ \\

If $(A,\mathfrak{m})$ is an Artinian local ring then $\Sp A$ with the operation $\mathfrak{m}^{i}\ast\mathfrak{m}^{j}=\mathfrak{m}^{i+j}$ can be viewed as a cyclic group of order $n$ where $n$ is the least natural number such that $\mathfrak{m}^{n}=0$. In particular, if $p$ is a prime number and $n\geq1$ then $\Sp A=\{\mathfrak{m},...,\mathfrak{m}^{n}\}$ is a cyclic group of order $n$ where $A=\mathbb{Z}/p^{n}\mathbb{Z}$ and $\mathfrak{m}=p\mathbb{Z}/p^{n}\mathbb{Z}$. \\

\end{document}